\documentclass[11pt]{amsart}
\usepackage[ansinew]{inputenc}
\usepackage{amssymb}
\usepackage{latexsym}
\usepackage{graphics}
\usepackage[dvips]{graphicx}
\usepackage{fancybox}
\usepackage
{amsmath,amsfonts}
\newtheorem{theorem}{Theorem}[section]
\newtheorem{lemma}[theorem]{Lemma}

\newtheorem{remark}{Remark}[section]
\newtheorem{ass}{Assumption}[section]

\begin{document}
\author{ M. M. Cavalcanti }
\address{ Department of Mathematics, State University of
Maring\'a, 87020-900, Maring\'a, PR, Brazil.}
\thanks{Research of Marcelo M. Cavalcanti partially supported by the CNPq Grant
300631/2003-0}

\author{ V. N. Domingos Cavalcanti }
\thanks{Research of Valéria N. Domingos Cavalcanti partially supported by the CNPq Grant
304895/2003-2}

\author{ R. Fukuoka }

\author{ J. A. Soriano}

\thanks{
distributed damping.}

\title[wave on compact manifolds]
{Asymptotic stability of the wave equation on compact surfaces and
locally distributed damping - A sharp result}

\maketitle
\begin{abstract}
This paper is concerned with the study of the wave equation on
compact surfaces and locally distributed damping, described by
\begin{equation}
\left.
\begin{array}{l}
u_{tt} -
\Delta_{\mathcal{M}}u+ a(x)\,g(u_{t})=0\, \ \;\,\,\,\,\text{%
on\ \thinspace }\mathcal{M}\times \left] 0,\infty \right[
,\smallskip
\end{array}%
\right.  \nonumber
\end{equation}
where $\mathcal{M}\subset \mathbb{R}^3$ is a smooth oriented
embedded compact surface without boundary. Denoting by
$\mathbf{g}$ the Riemannian metric induced on $\mathcal{M}$ by
$\mathbb{R}^3$, we prove that for each $\epsilon
> 0$, there exist an open subset $V \subset\mathcal M$
and a smooth function $f:\mathcal M \rightarrow \mathbb R$ such
that $meas(V)\geq meas(\mathcal M)-\epsilon$, $Hess f \approx
\mathbf{g}$ on $V$ and $\underset{x\in V}\inf |\nabla f(x)|>0$.

In addition, we prove that if $a(x) \geq a_0> 0$ on an open subset
$\mathcal{M}{\ast} \subset \mathcal M$ which contains
$\mathcal{M}\backslash V$ and if $g$ is a monotonic increasing
function such that $k |s| \leq |g(s)| \leq K |s|$ for all $|s|
\geq 1$, then uniform and optimal decay rates of the energy hold.
\end{abstract}

\section{Introduction}

\setcounter{equation}{0}

Let $\mathcal{M}$ be a smooth oriented embedded compact surface
without boundary in $\mathbb{R}^{3}$ and let $\mathbf g$ denote
the Riemannian metric induced on $\mathcal M$ by $\mathbb R^3$.
For $\epsilon> 0$ we prove that there exist an open subset $V
\subset\mathcal M$ and a smooth function $f:\mathcal M \rightarrow
\mathbb R$ such that $meas(V)\geq meas(\mathcal M)-\epsilon$,
$Hess f \approx \mathbf{g}$ on $V$ and $\underset{x\in V}\inf
|\nabla f(x)|>0$ (See Subsection
\ref{hessianofuncaoquasemetrica}).

We denote by $\nabla_{T}$ the tangential-gradient on $\mathcal{M}
$ and by $\Delta_{\mathcal{M}}$ the Laplace-Beltrami operator on
$\mathcal{M}$. This paper is devoted to the study of the uniform
stabilization of solutions of the
following damped problem%
\begin{equation}
\left\{
\begin{array}{l}
u_{tt}-\Delta_{\mathcal{M}}u+ a(x)\,g(u_{t})=0\,\,\,\ \ \ \ \ \ \
\ \ \ \ \ \ \ \ \ \ \ \ \;\,\,\,\,\text{
on\ \thinspace }\mathcal{M} \times \left] 0,\infty \right[ ,\smallskip \\
u(x,0)=u^0(x),\quad u_t(x,0)=u^1(x)\ \ \ \ \ \ \ \ \ \ \ \ \ \ \ \
x\in \mathcal{M},
\end{array}%
\right.  \label{1.1}
\end{equation}
where $a(x) \geq a_0 > 0$ on an open proper subset
$\mathcal{M}_{\ast} \supset \mathcal M \backslash V$ of
$\mathcal{M}$ and in addition $g$ is a monotonic increasing function
such that $k |s| \leq |g(s)| \leq K |s|$ for all $|s| \geq 1$.

\smallskip

A natural question arises in the context of the wave equation on
compact surfaces: {\em Would it be possible to stabilize the
system by considering a localized feedback acting only on a
portion of the surface ?} In affirmative case, {\em what would be
the geometrical impositions we have to assume on the surface?}
When the damping term acts on the whole surface, the conjecture
was studied by Cavalcanti and Domingos Cavalcanti in
\cite{Cavalcanti} and also by Andrade et al. in \cite{Cavalcanti3,
Cavalcanti4} in the context of viscoelastic problems. For linear
waves, we can mention the works due to Rauch and
Taylor\cite{Rauch}, Hitrik \cite{Hitrik} and, recently
Christianson \cite{Christianson}. For the nonlinear wave equation
on compact manifolds with boundary, it is important to cite the
work due to Triggiani and Yao \cite{Trigianni}. More recently, the
authors of the present work \cite{Cavalcanti6} studied the linear
wave equation on a compact surface $\mathcal{M}$ without boundary
supplemented with a nonlinear and localized dissipation. In this
article the authors prove the above conjecture when the portion of
$\mathcal{M}$ where the damping is effective is strategically
chosen. Setting $\mathcal{M}=\mathcal{M}_0 \cup \mathcal{M}_1 $,
where
\begin{eqnarray*}
\mathcal{M}_1:=\{x\in \mathcal{M}; m(x) \cdot \nu(x) > 0\} ~\hbox{
\and }\mathcal{M}_0=\mathcal{M}\backslash \mathcal{M}_1,
\end{eqnarray*}
$m(x):=x-x^0$, ($x^0 \in \mathbb{R}^3$ fixed) and $\nu$ is the
exterior unit normal vector field of $\mathcal{M}$, then for
$i=1,\ldots,k$, they assume that there exist open subsets
$\mathcal M_{0i} \subset \mathcal M_0$ of $\mathcal M$ with smooth
boundary $\partial \mathcal M_{0i}$ such that $\mathcal M_{0i}$
are umbilical. Moreover, they suppose that the {\em mean
curvature} $H$ of each $\mathcal{M}_{0i}$ is {\em non-positive}
(i.e. $H\leq 0$ on $\mathcal{M}_{0i}$ for every $i=1,\ldots,k$)
and that the damping is effective on an open subset
$\mathcal{M}_{\ast} \subset \mathcal M$ which contains
$\mathcal{M}\backslash\cup_{i=1}^k \mathcal M_{0i}$. Roughly
speaking, the region which does not contain dissipative effects
must be {\em umbilical}. This is required since the authors employ
the same multipliers considered to solve the similar question for
the wave equation,
\begin{eqnarray*}
u_{tt}- \Delta u + a(x) g(u_t) =0 ~ \hbox{ in }~\Omega \times
(0,+\infty),
\end{eqnarray*}
where $\Omega$ is a bounded domain of $\mathbb{R}^n$ with smooth
boundary $\partial \Omega$. They considered the well known
multiplier given by the vector field $m(x):= x-x^0$, $x^0\in
\mathbb{R}^n$ arbitrarily chosen, but fixed, taken out of the
domain $\Omega$, according to the figure 1 below.

%%%%%%%%%%%%%%%%%%%%%%%%%%%%%%%%%%%%%%% BEGINNING OF FIGURE 1 %%%%%%%%%%%%%%%%%%%%%%
{\small
\begin{figure}[ht]
\begin{center}
\setlength{\unitlength}{0.73pt}
\begin{picture}(440,210)(0,53)

\put(40,170){$\mathcal{M}_1$} \put(350,155){$\mathcal{M}_{0}$}
%\put(90,160){\LARGE$\omega$}

%\put(200,160){\LARGE$\Omega\backslash\omega$}

\put (460,160){$x_0$} \put(295,234){$\circ$} \put(295,85){$\circ$}
\put(310,237){$\mathcal{M}_{\ast}$}
\put(440,167){\line(-2,1){140}} \put(440,167){\line(2,-1){14}}
\put(455,166){\line(-2,-1){150}} \put(440,165){\line(-5,1){380}}
\put(67,239){\vector(-3,1){12}} \put(319,189){\vector(-4,1){12}}
\put(297,198){$x-x^0$}

\put(344,183){\vector(2,1){25}} \put(350,198){$\nu(x)$}
\put(106,232){\vector(-1,1){25}} \put(20,245){$x-x^0$}
\put(90,255){$\nu(x)$} \put(440,165){\line(5,-1){15}}
\put(445,160){$\bullet$}

\bezier{540}(347,163)(347,252)(207,253)
\bezier{740}(67,163)(67,253)(207,253)
\bezier{740}(347,163)(347,73)(207,73)
\bezier{740}(67,163)(67,73)(207,74)
\put(60,43){\setlength{\unitlength}{0.59pt}
%%%%%%%%%%%%%%%%%%%%%%%%%%%%%%
\bezier{740}(9,150)(20,160)(30,165)
%\bezier{740}(30,165)(35,168)(50,171)
\bezier{740}(50,171)(55,172)(70,173)
%\bezier{740}(70,173)(80,173)(90,174)
\bezier{740}(90,174)(100,175)(115,175)
%\bezier{740}(115,175)(125,176)(140,176)
\bezier{740}(140,176)(155,177)(165,177)
%\bezier{740}(165,177)(175,177)(185,177)
\bezier{740}(185,177)(200,177)(210,176)
%\bezier{740}(210,176)(225,176)(235,175)
\bezier{740}(235,175)(250,175)(255,174)
%\bezier{740}(255,174)(270,174)(280,173)
\bezier{740}(280,173)(295,173)(305,172)
%\bezier{740}(305,172)(320,172)(330,170)
\bezier{740}(330,170)(345,168)(345,167)
\bezier{740}(345,167)(350,166)(354,163)
%%%%%%%%%%%%%%%%%%%%%%%%%%%%%%%%%%%%%%%%
\bezier{740}(9,150)(20,130)(175,130)
\bezier{740}(175,130)(330,130)(354,163)
%%%%%%%%%%%%%%%%%%%%%%%%%%%%%%%%%%%%
%\bezier{740}(295,234)

\bezier{740}(295,56)(300,60)(305,72)
%\bezier{740}(305,72)(307,76)(312,88)
\bezier{740}(312,88)(313,92)(317,102)
%\bezier{740}(317,102)(319,106)(321,116)
\bezier{740}(321,116)(322,120)(324,130)
%\bezier{740}(324,130)(325,134)(327,144)
\bezier{740}(327,144)(328,148)(329,155)

\bezier{740}(329,155)(329,162)(327,169)
%\bezier{740}(327,169)(326,174)(324,180)
\bezier{740}(324,181)(323,184)(315,194)
%\bezier{740}(315,196)(313,200)(308,210)
\bezier{740}(308,210)(306,214)(303,220)
%\bezier{740}(303,220)(301,224)(298,230)
\bezier{740}(298,230)(296,234)(295,238)
%%%%%%%%%%%%%%%%%%%%%%%%
\bezier{740}(295,56)(250,150)(295,238)
 }

\put(180,40){Figure 1}
\end{picture}
\end{center}
\caption{The observer is at $x_0$. The subset $\mathcal M_0$ is
the ``visible'' part of $\mathcal M$ and $\mathcal M_1$ is its
complement. The subset $\mathcal{M}_{\ast}\supset \mathcal
M\backslash \cup_{i=1}^k \mathcal M_{0i}$ is an open set which
contains $\mathcal M\backslash \cup_{i=1}^k \mathcal M_{0i}$ and
the damping is effective there. Observe that in figure 1, $k=1$
and $\mathcal{M}_{0i}=\mathcal{M}_{01}=\mathcal{M}_0$.}

\end{figure}

}
%%%%%%%%%%%%%%%%%%%%%%%%%%%% END OF FIGURE %%%%%%%%%%%%%%%%%%%%%%%%%%%%%%%%%%%

\smallskip

Once the multiplier $m(x) = x- x^0$ is not {\em intrinsically
connected } with the manifold $\mathcal{M}$ they have to impose a
restriction on the part $\mathcal{M}_0$ (without damping), namely,
$\mathcal{M}_0$ must be {\em umbilical, or umbilical by parts}.

\smallskip

The main goal of the present manuscript is to improve considerably
the previous result due to \cite{Cavalcanti6}, {\em reducing
arbitrarily} the volume of the region where the dissipative effect
lies. For this purpose we will construct an intrinsic multiplier
that will play a crucial role when establishing the desired
uniform decay rates of the energy. Fix $\epsilon > 0$. This
multiplier is, roughly speaking, given by the $\nabla_T f$, where
$f:\mathcal{M} \rightarrow \mathbb{R}$ is a regular function which
verifies $Hess f \approx \mathbf g$ and $\underset{x\in
V}\inf|\nabla f (x)|>0$ on a subset $V$ of $\mathcal{M}$ such that
$meas(V)\geq meas (\mathcal{M})-\epsilon$. This construction will
be clarified in subsections 4.3 and 4.4.

\medskip

We would like to emphasize that the proofs of \cite{Rauch,
Bardos,Hitrik}, based on microlocal analysis, {\em do not extend}
to the nonlinear problem (\ref{1.1}). In addition, making use of
arguments due to Cavalcanti, Domingos Cavalcanti and Lasiecka
\cite{Cavalcanti5}, we obtain {\it explicit and optimal decay
rates of the energy}. The obtained decay rates are optimal, since
they are the same as these optimal rates derived in the works of
Alabau-Boussouira \cite{Alabau} or Toudykov \cite{Daniel}.

Our paper is organized as follows.   Section 2 is concerned with
the statement of the problem  and we introduce some notation . Our
main result is stated in Section 3. Section 4 is devoted to the
proof of the main result.

\section{Statement of the Problem}

\setcounter{equation}{0}

Let $\mathcal{M}$ be a smooth oriented embedded compact surface
without boundary in $\mathbb{R}^{3}$. For $\epsilon >0$ we prove
that there exist an open subset $V \subset\mathcal M$ and a smooth
function $f:\mathcal M \rightarrow \mathbb R$ such that
$meas(V)\geq meas(\mathcal M)-\epsilon$, $Hess f \approx
\mathbf{g}$ on $V$ and $\underset{x\in V}\inf|\nabla f(x)|>0$ (See
Subsection \ref{hessianofuncaoquasemetrica}).

In this paper, we investigate the stability properties of function
$u(x,t)$ which solves the damped problem
\begin{equation}
\left\{
\begin{array}{l}
u_{tt}- \Delta_{\mathcal{M}}u+ a(x)\,g(u_{t})=0\, \ \ \ \text{
on\ \thinspace }\mathcal{M}\times \left] 0,\infty \right[ ,\smallskip \\
u(0)=u^0, \quad u_t(0)=u^1,\ \ \ \ \ \
\end{array}%
\right.  \label{3.1}
\end{equation}
where the feedback function $g$ \ satisfies the Assumption 2.1.

\begin{ass}\label{ass:1}
${}$
\smallskip

$(i) $ $\ \ \ g\left( s\right) $ {\em is continuous and monotone
increasing},

$( ii) $ $\ \ g\left( s\right) s>0$ {\em for} $%
s\neq 0,$

$\left( iii\right) $ $\ \ k\,|s|\leq g\left( s\right) \leq K\,|s
|$ {\em for} $\left\vert s\right\vert >1,$
\end{ass}
\noindent where $k$ and $K$ are two positive constants.

In addition, to obtain the stabilization of problem $(\ref{3.1}),$
we shall need the following geometrical assumption:
\begin{ass}\label{as:2}
Assume that $a \in L^{\infty}(\mathcal{M})$ is a nonnegative
function such that
\begin{eqnarray}\label{eq:2.2}
a(x) \geq a_0 > 0, \quad \hbox{\em a. e. on  }\mathcal{M}_{\ast},
\end{eqnarray}
where $\mathcal{M}_{\ast}$ is an open set of $\mathcal{M}$ which
contains $\mathcal{M}\backslash V$.

\end{ass}

\medskip
In the sequel, we are going to consider $\Sigma =\mathcal{M}
\times \left] 0,T\right[ $ and the Sobolev spaces
$H^s(\mathcal{M})$,\,$s\in \mathbb{R}$, as in Lions and Magenes
\cite{Magenes} section 7.3.

On the other hand, using the Laplace-Beltrami operator
$\Delta_{\mathcal{M}}$ on $\mathcal{M}$, we can give a more
intrinsic definition of the spaces $H^s(\mathcal{M})$. Considering
\begin{eqnarray*}
H^{2m}\left( \mathcal{M} \right) &=&\left\{ u\in L^2(\mathcal{M}
)\,/\Delta_{\mathcal{M}}^m \,u \in L^2(\mathcal{M})\right\} ,
\end{eqnarray*}
which, equipped with the canonical norm
\begin{equation}\label{norm H^s II}
\left\Vert u\right\Vert_{H^{2m}(\mathcal{M})}^{2}=\left\Vert
u\right\Vert _{L^{2}(\mathcal{M} )}^{2}+\left\Vert
\Delta_{\mathcal{M}}^m u\right\Vert _{L^{2}(\mathcal{M} )}^{2},
\end{equation}
is a Hilbert space.
\smallskip

We set
\begin{eqnarray*}
&V:= \{v\in H^1(\mathcal{M}); \int_{\mathcal{M}}
v(x)\,d\mathcal{M} =0 \},&
\end{eqnarray*}
which is a Hilbert space with the topology endowed by
$H^1(\mathcal{M})$.

Note that the condition $\int_{\mathcal{M}} v(x)\,d\mathcal{M} =0$
is required in order to guarantee the validity of the Poincaré
inequality,
\begin{eqnarray}\label{Poincare}
||f||_{L^2(\mathcal{M})}^2 \leq (\lambda_1)^{-1} ||\nabla_T
f||_{L^2(\mathcal{M})}^2, \quad \hbox{ for all } f\in V,
\end{eqnarray}
where $\lambda_1$ is the first eigenvalue of the Laplace-Beltrami
operator.

We observe that problem (\ref{3.1}) can be rewritten as
$$
\frac{dU}{dt} + \mathcal{A}U = G(U),
$$
where
\begin{equation*}
\mathcal{A}=\left(
\begin{array}{c}
\,\,0\,\,\,\,\,\,\,\,-I \\
-\Delta_{\mathcal{M}}\,\,\,\,\,\,\,\quad 0
\end{array}
\right)
\end{equation*}
is a maximal monotone operator and $G(\cdot)$ represents a locally
Lipschitz perturbation. So, making use of standard semigroup
arguments we have the following result:

\medskip

\begin{theorem}
${}$
\begin{itemize}
\item $\left( \mathbf{i}\right) $ \textit{Under
the above-mentioned conditions, problem }$\left( \ref{3.1}\right) $\textit{\ is wellposed in%
}$\,\,$\textit{the space}$\,\,V\times L^2(\mathcal{M})$\textit{,\
that is, for any initial data \thinspace }$\left\{
u^{0},u^{1}\right\} \in V\times L^2(\mathcal{M})$\textit{, there
exists a unique\thinspace weak solution
of }(\ref{3.1})\textit{\ in the class}%
\begin{equation}
u \in C(\mathbb{R}_{+};V)\cap C^{1}(\mathbb{R}_{+};%
L^2(\mathcal{M})).  \label{3.9}
\end{equation}

\item $ \left(\mathbf{ii}\right) $\textit{In addition, the
velocity of the solution has the regularity given by}
\begin{equation}
u_{t}\in  L_{loc}^{2}\left( \mathbb{R}_{+};L^{2}\left( \mathcal{M}
\right) \right) , \label{3.10}
\end{equation}%
\textit{and, consequently,} $g\left( u_{t}\right) \in
L_{loc}^{2}\left( \mathbb{R}_{+};L^{2}\left( \mathcal{M}\right)
\right) $ by Assumption \ref{ass:1}.
\end{itemize}
\noindent \textit{Furthermore, if} $\left\{ u^{0},u^{1}\right\}
\in  V\cap H^{2}\left( \mathcal{M} \right) \times V $
\textit{then, the solution has the following regularity:}
\begin{eqnarray*}
 u \in
L^{\infty }\left( \mathbb{R}_{+};V\cap H^{2}\left( \mathcal{M}
\right) \right) \cap W^{1,\infty }\left(\mathbb{R}_{+}; V
\right)\cap W^{2,\infty }\left( \mathbb{R}_{+};L^{2}\left(
\mathcal{M} \right) \right) .
\end{eqnarray*}
\end{theorem}
\smallskip

Consider that $u $ is the unique global weak solution of problem
(\ref{3.1}) given by Theorem 2.1. We define the corresponding
energy functional by
\begin{equation}
E(t)=\frac{1}{2}\int_{\mathcal{M}}\left[ \left\vert
u_{t}(x,t)\right\vert ^{2}+\left\vert \nabla _T  u(x,t)\right\vert
^{2}\right] d\mathcal{M} . \label{3.11}
\end{equation}

For every solution of (\ref{3.1}), in the class (\ref{3.9}) we
obtain for all $t_2 > t_1 \geq 0$
\begin{equation}
E(t_{2})-E(t_{1})=-\int_{t_{1}}^{t_{2}}\int_{\mathcal{M}
}a(x)\,g(u_{t})u_{t}\,d\mathcal{M} dt, \label{3.12}
\end{equation}%
and, therefore, the energy is a non increasing function of the
time variable $ t $.

\section{Main Result}

\setcounter{equation}{0}

In order to state the stability result, we need to define some
functions which are firstly introduced in Lasiecka and Tataru
\cite{Lasiecka-Tataru}. For the reader's comprehension we will
repeat them briefly.
Let $h$ be a\ concave, strictly increasing function, with $%
h\left( 0\right) =0$, and such that
\begin{equation}
h\left( s\,g(s))\right) \geq s^{2}+g^{2}(s),\text{\ for
}\left\vert s\right\vert \leq 1.  \label{4.1'}
\end{equation}

Note that such function can be straightforwardly constructed,
considering the hypotheses on $g$ in Assumption \ref{ass:1}. In
view of this function, we define
\begin{equation}
r(.)=h(\frac{.}{meas\left( \Sigma _{1}\right) }).  \label{4.2'}
\end{equation}
As $r$ is monotone increasing, then $cI+r$ is invertible for all
$c\geq 0.$ For $L$ a positive constant, we set
\begin{equation}
p(x)=(cI+r)^{-1}\left( Lx\right) ,  \label{4.3'}
\end{equation}%
where the function $p$ is easily seen to be positive, continuous
and strictly increasing with $p(0)=0$. Finally, let
\begin{equation}
q(x)=x-(I+p)^{-1}\left( x\right) .  \label{4.4'}
\end{equation}%
We can now proceed to state our stability result.

\medskip

\begin{theorem}\label{Theo. 3.1}
Assume that Assumption \ref{ass:1} and Assumption \ref{as:2} are
in place. Let $u$ be the weak solution of problem (\ref{3.1}).
With the energy $E(t)$ defined as in (\ref{3.11}), there exists
$T_{0}>0$ such that
\begin{equation}
E(t)\leq S\left( \frac{t}{T_{0}}-1\right) ,\text{ \ }\forall
t>T_{0}, \label{4.5'}
\end{equation}%
with $\underset{t\rightarrow \infty }{\lim }S(t)=0,$ where the
contraction semigroup $S(t)$ is the solution of the differential
equation
\begin{equation}
\frac{d}{dt}S(t)+q(S(t))=0,\text{\ \ }S(0)=E(0),  \label{4.6'}
\end{equation}
where $q$ is given in (\ref{4.4'}), the constant $L$, which is
given in (\ref{4.3'}), depends on $meas(\Sigma)$ and the constant
$c$ is equal to $\frac{ k^{-1}+K}{meas\left( \Sigma\right)
(1+||a||_{\infty})}.$
\end{theorem}
\smallskip

\noindent \begin{remark} If the feedback is linear, e. g., $
g(s)=s,$ then, under the same assumptions of Theorem \ref{Theo.
3.1}, we obtain that the energy of problem (\ref{3.1}) decays
exponentially with respect to the initial energy, that is, there
exist two positive constants $C>0$ and $k>0$ such
that%
\begin{equation}
E(t)\leq Ce^{-kt}E(0),\text{ \ \ }t>0.  \label{4.7'}
\end{equation}

\smallskip

If we consider $g(s) = s^ p $, $ p> 1 $ at the origin, and since
the function $s^{\frac{p+1}{2}} $ is convex for $p > 1 $, then,
solving
\begin{equation}\label{S1}
S_t + S^{\frac{p+1}{2}} =0
\end{equation}
we obtain the following polynomial decay rate:
$$ E(t) \leq C(E(0))[E(0)^{\frac{-p +1}{2}} +
t(p-1)]^{\frac{ 2}{-p +1}},$$

We can find more examples of explicit decay rates in Cavalcanti,
Domingos Cavalcanti and Lasiecka \cite{Cavalcanti5}.
\end{remark}

\section{Proof of Main result}

\setcounter{equation}{0}

\subsection{Preliminaries}
${}$
\smallskip

\smallskip

We collect, below, some few formulas to be invoked in the sequel.

\smallskip

Let $\nu$ be the exterior normal vector field on $\mathcal{M}$.
The {\em Laplace- Beltrami operator} $\Delta _{\mathcal{M}}$ of a
function $\varphi :\mathcal{M}\rightarrow \mathbb{R}$ of class
$C^{2}$ is defined by
\begin{equation}
\Delta _{\mathcal{M}}\varphi :=div_{T}\nabla _{T}\varphi ,
\label{4.4}
\end{equation}
where $div_{T}\nabla_{T}\varphi$, is the $divergent$ of the vector
field $\nabla_{T}\varphi$.

Assuming that $\varphi :\mathcal{M}\rightarrow \mathbb{R}$ is a
function of class $C^{1}$ and $x\in \mathcal{M} \mapsto q(x) \in
T_x(\mathcal{M})$ is a vector field of class $C^{1}$, we have,
\begin{eqnarray}
\int_{\mathcal{M} }q\cdot \nabla _{T}\varphi \,d\mathcal{M}
&=&-\int_{\mathcal{M} }div_T q\,\varphi \,d\mathcal{M},
\label{4.5}
\end{eqnarray}
\begin{eqnarray}
2\varphi (q\cdot \nabla_T\varphi)= q \cdot
\nabla_T(\varphi^2).\label{E}
\end{eqnarray}

From (\ref{4.5}) and (\ref{E}), we conclude the following formula
\begin{equation}
2\int_{\mathcal{M} }\varphi (q\cdot \nabla _{T}\varphi
)\,d\mathcal{M} =\int_{\mathcal{M} }q\cdot \nabla _{T}(\varphi
^{2})\,d\mathcal{M} =-\int_{\mathcal{M} }div_{T}\,q|\varphi
|^{2}d\mathcal{M} . \label{second formula}
\end{equation}

\smallskip

We define a continuous linear operator $-\Delta
_{\tilde{\mathcal{M}}}:H^{1}(\tilde{\mathcal{M}})\rightarrow
(H^{1}(\tilde{\mathcal{M}}))^{\prime }$, where
$\tilde{\mathcal{M}}$ is a nonempty open subset of $\mathcal{M}$
(sometimes the whole $\mathcal{M}$) such that

\begin{equation}
\langle -\Delta _{\tilde{\mathcal{M}}}\varphi ,\psi \rangle
=\int_{\tilde{\mathcal{M}}}\nabla _{T}\varphi \cdot \nabla
_{T}\psi \,d\mathcal{M} , \quad \forall \varphi, \psi \in
H^{1}(\tilde{\mathcal{M}}) \label{4.12'}
\end{equation}%
and, in particular,
\begin{equation}
\langle -\Delta _{\tilde{\mathcal{M}}}\varphi ,\varphi \rangle
=\int_{\tilde{\mathcal{M}}}|\nabla _{T}\varphi
|^{2}\,d\mathcal{M}, \quad \forall \varphi \in
H^{1}(\tilde{\mathcal{M}}). \label{4.13'}
\end{equation}

\smallskip

The operator $-\Delta_{\tilde{\mathcal{M}}}+I$ defines an
isomorphism from $H^1(\tilde{\mathcal{M}})$ over
$[H^1(\tilde{\mathcal{M}})]'$. We observe that when
$\tilde{\mathcal{M}}$ is a manifold without boundary, and this is
the case, for instance, if $\tilde{\mathcal{M}}=\mathcal{M}$, we
have $H^1(\tilde{\mathcal{M}})= H_0^1(\tilde{\mathcal{M}})$ and,
consequently,
$[H^1(\tilde{\mathcal{M}})]'=H^{-1}(\tilde{\mathcal{M}})$.

\begin{remark}
It is convenient to observe that all the above classical formulas
can be extended to Sobolev spaces using density arguments.
\end{remark}

The proof of Theorem 3.1 proceeds through several steps. In order
to obtain the decay rate stated in (\ref{4.5'}), we will consider,
initially, regular solutions of problem (\ref{3.1}). Then, making
use of standard density arguments, the estimate (\ref{4.5'}) holds
for weak solutions.

\subsection{An identity}

We begin proving the following proposition:

\smallskip

\noindent \textbf{Proposition 4.2.1. }\textit{Let}
$\mathcal{M}\subset \mathbb{R}^3 $ \textit{be an oriented regular
compact surface without boundary} \textit{\ and }$q$ \textit{a
vector field of class } $C^1.$ \textit{\ Then, for every regular
solution }$u $\textit{\thinspace of }(\ref{1.1})\textit{ we have
the following identity}
\begin{eqnarray}
&&\left[ \int_{\mathcal{M}}u_{t}\, q\cdot \nabla_T
u\,d\mathcal{M}\right]
_{0}^{T}+\frac{1}{2}\int_{0}^{T}\int_{\mathcal{M}
}(div_{T}q)\left\{ \left\vert u_{t}\right\vert ^{2}-\left\vert
\nabla _{T}
u\right\vert^{2}\right\} d\mathcal{M} dt \label{5.2.1} \\
&&+\int_{0}^{T}\int_{\mathcal{M}} \nabla_T u \cdot \nabla_T q\cdot
\nabla_T u \,d\mathcal{M} dt +\int_{0}^{T}\int_{\mathcal{M}
}a(x)\,g(u_{t})(q\cdot\nabla_T u)d\mathcal{M} dt=0. \nonumber
\end{eqnarray}
\noindent \textbf{Proof. } Multiplying the  equation (\ref{1.1})
by the multiplier $q \cdot \nabla_T u$ and integrating on
$\mathcal{M}\times ]0,T[$, we obtain
\begin{eqnarray}\label{IDEN}
0= \int_0^T \int_{\mathcal{M}} (u_{tt} - \Delta_{\mathcal{M}} u +
a(x)g(u_t))(q \cdot \nabla_T u)\,d\mathcal{M}\,dt.
\end{eqnarray}

Next, we will estimate some terms on the RHS of identity
(\ref{IDEN}).  Taking (\ref{4.5}), (\ref{E}) and (\ref{second
formula}) into account, we obtain
\begin{eqnarray}
&&\int_{0}^{T}\int_{\mathcal{M}}\left( -\Delta_{\mathcal{M}}
u\right) (q\cdot \nabla_T u)\,d\mathcal{M} dt
=\int_{0}^{T}\int_{\mathcal{M}} \nabla_{T} u \cdot\nabla
_{T}(q\cdot \nabla_T u)\, d\mathcal{M}
dt   \label{5.2.4} \\
&&=\int_{0}^{T}\int_{\mathcal{M}}\nabla _{T} u \cdot \nabla _{T} q
\cdot \nabla_{T} u \,d\mathcal{M} dt
+\frac{1}{2}\int_{0}^{T}\int_{\mathcal{M}}q \cdot \nabla_T
[|\nabla_T u|^2]d\mathcal{M}
dt  \notag \\
&&=\int_{0}^{T}\int_{\mathcal{M}}\nabla _{T} u \cdot \nabla _{T} q
\cdot \nabla_{T} u \,d\mathcal{M} dt
-\frac{1}{2}\int_{0}^{T}\int_{\mathcal{M}}\left\vert  \nabla _{T}
u\right\vert^{2}div_{T}q\,d\mathcal{M} dt,  \notag
\end{eqnarray}
and, integrating by parts and considering (\ref{second formula}),
we obtain
\begin{eqnarray}
&&\int_{0}^{T}\int_{\mathcal{M}}\left( u_{tt}+ a(x)\,g( u_{t})
\right)
(q \cdot \nabla_T u)\,d\mathcal{M} dt \label{5.2.5}\\
&=&\left[ \int_{\mathcal{M}} u_t(q \cdot \nabla_T u)\right]
_{0}^{T}- \int_{0}^{T}\int_{\mathcal{M}}u_t (q\cdot \nabla_T u_t)d\mathcal{M} dt  \notag \\
&&+\int_{0}^{T}\int_{\mathcal{M}}a(x)\,g\left(
u_{t}\right) (q \cdot \nabla_T u)d\mathcal{M} dt  \notag \\
&=&\left[ \int_{\mathcal{M}}u_{t} (q \cdot \nabla_T u)\right]
_{0}^{T}+\frac{1}{2} \int_{0}^{T}\int_{\mathcal{M}
}(div_{T}q_{T})\left\vert u_{t}\right\vert
^{2}d\mathcal{M} dt  \notag \\
&&+\int_{0}^{T}\int_{\mathcal{M}}a(x)\,g\left( u_{t}\right) (q
\cdot \nabla_T u)d\mathcal{M} dt. \notag
\end{eqnarray}

Combining (\ref{IDEN}), (\ref{5.2.4}) and ( \ref{5.2.5}), we deduce
(\ref{5.2.1}), which concludes the proof of Proposition 4.2.1.
$\quad \square $

\medskip

Employing (\ref{5.2.1}) with $q(x)= \nabla_T f$ where
$f:\mathcal{M} \rightarrow \mathbb{R}$ is a $C^3$ function to be
determined later, we infer
\begin{eqnarray}
&&\left[ \int_{\mathcal{M}}u_{t}\,\nabla_T f \cdot \nabla_T
u\,d\mathcal{M}\right] _{0}^{T} + \frac12
\int_{0}^{T}\int_{\mathcal{M} }\Delta_{\mathcal{M}}f\left\{
\left\vert u_{t}\right\vert ^{2}-\left\vert
\nabla _{T} u\right\vert^{2}\right\} d\mathcal{M} dt \label{first inequality}\\
&&+\int_{0}^{T}\int_{\mathcal{M}} (\nabla_T u \cdot Hess(f) \cdot \nabla_T u)\,d\mathcal{M} dt  \notag \\
&&+\int_{0}^{T}\int_{\mathcal{M}}a(x)\,g(u_{t})(\nabla_{T}
f\cdot\nabla_T u)d\mathcal{M} dt=0. \nonumber
\end{eqnarray}

\smallskip

\noindent We have the following identity:

\smallskip

\noindent \textbf{Lemma 4.2.3. }\textit{Let }$u $\textit{\ be a
weak solution to problem }(\ref{1.1}) \textit{and} $\xi\in
C^1(\mathcal{M})$.\textit{\ Then}
\begin{eqnarray}\qquad
\left[\int_{\mathcal{M}}u_t\,\xi\, u\,d\mathcal{M}  \right]_0^T
&=& \int_0^T\int_{\mathcal{M}} \xi |u_t|^2 d\mathcal{M} dt-
\int_0^T \int_{\mathcal{M}} \xi |\nabla_T u|^2 d\mathcal{M} dt \label{identity II}\\
&-&  \int_0^T \int_{\mathcal{M}} (\nabla_T u \cdot \nabla_T
\xi)u\,d\mathcal{M} dt -\int_0^T \int_{\mathcal{M}}
a(x)\,g(u_t)\,\xi\,u\,d\mathcal{M} dt.\nonumber
\end{eqnarray}
\textbf{Proof:} Multiplying the equation of (\ref{1.1}) by $\xi
\,u$ and integrating by parts we obtain the desired result. \quad
$\square$
\smallskip

Substituting $\xi=\alpha>0$ in (\ref{identity II}) and combining
the obtained result with identity (\ref{first inequality}) we
deduce
\begin{eqnarray}
&&\int_{0}^{T}\int_{\mathcal{M}}(\frac{\Delta_{\mathcal{M}}f}{2}-\alpha)
\left\vert u_{t}\right\vert
^{2} d\mathcal{M} dt.\notag\\
&&+\int_0^T \int_{\mathcal{M}} \left[  (\nabla_T u \cdot Hess(f)
\cdot \nabla_T
u)+\left(\alpha-\frac{\Delta_{\mathcal{M}}f}{2}\right)|\nabla_T
u|^2\right] \,d\mathcal{M} dt\nonumber\\
&&=-\left[ \int_{\mathcal{M}}u_{t}\, \nabla_{T} f\cdot \nabla_T
u\,d\mathcal{M}\right] _{0}^{T}-\alpha\left[
\int_{\mathcal{M}}u_{t}\, u\,d\mathcal{M}\right] _{0}^{T} \label{second inequality} \\
&& -\alpha\int_0^T \int_{\mathcal{M}}a(x)\,g(u_t)u\, d\mathcal{M}
dt\nonumber\\
&&-\int_{0}^{T}\int_{\mathcal{M} }a(x)\,g(u_{t})(\nabla_{T}
f\cdot\nabla_T u)d\mathcal{M} dt. \nonumber
\end{eqnarray}

\begin{remark}
This is the precise moment where the properties of function $f$
play an important role. Note that what we just need is to find a
subset $V$ of $\mathcal{M}$ such that
\begin{eqnarray}\label{main inequality}
&& C \int_0^T \int_{V}\left[ u_t^2 + |\nabla_T u|^2 \right]
d\mathcal{M}dt\\
&&\leq\int_{0}^{T}\int_{V}(\frac{\Delta_{\mathcal{M}}f}{2}-\alpha)
\left\vert u_{t}\right\vert
^{2} d\mathcal{M} dt\nonumber\\
&&+\int_0^T \int_{V} \left[  (\nabla_T u \cdot Hess(f) \cdot
\nabla_T
u)+\left(\alpha-\frac{\Delta_{\mathcal{M}}f}{2}\right)|\nabla_T
u|^2\right] \,d\mathcal{M} dt,\nonumber
\end{eqnarray}
for some positive constant $C$, provided that $\alpha$ is suitably
chosen. Assuming, for a moment, that (\ref{main inequality})
holds,  (\ref{second inequality}) yields
\begin{eqnarray}\label{B}
2C \int_{0}^{T} E(t) dt &\leq& C \int_0^T
\int_{\mathcal{M}\backslash V}\left[ u_t^2 +
|\nabla_T u|^2 \right]d\mathcal{M}dt \nonumber\\
&+& \left|\left[ \int_{\mathcal{M}}u_{t}\, \nabla_{T} f\cdot
\nabla_T u\,d\mathcal{M}\right] _{0}^{T}\right|+\alpha\left|\left[
\int_{\mathcal{M}}u_{t}\, u\,d\mathcal{M}\right] _{0}^{T}\right| \label{A} \\
&+& \left|\alpha\int_0^T \int_{\mathcal{M}}a(x)\,g(u_t)u\,
d\mathcal{M}dt\right|\nonumber\\
&+&\left|\int_{0}^{T}\int_{\mathcal{M} }a(x)\,g(u_{t})(\nabla_{T}
f\cdot\nabla_T u)d\mathcal{M} dt\right|. \nonumber
\end{eqnarray}

The inequality (\ref{A}) is controlled considering a standard
procedure, which, for the reader's convenience, we will repeat
later. The main idea behind this procedure is to consider the
dissipative area, namely, $\mathcal{M}_{\ast}$, containing the set
$\mathcal{M}\backslash V$ as stated in (\ref{eq:2.2}). It is
important to observe that $\mathcal{M}_{\ast}$ is as small as big
$V$ can be.
\end{remark}

\medskip

The next subsections are devoted to the construction of a function
$f$ as well as a subset  $V$ of $\mathcal{M}$ such that the
inequality (\ref{main inequality}) holds. This will be done, for
simplicity, in a general setting, that is, for a Riemannian
manifold (without boundary) with Riemmanian metric $\mathbf{g}$ of
class $C^2$.

%%%%%%%%%%%%%%%%%%%%%%%%%%%%%%%%%%%%%% Hessian almost diagonal %%%%%%%%%%%%%%%%%%%%%%%%%%%%%%%%%%%%%%%%%%%%%%%%%%%%%%%%%

\subsection{Construction of a function such that $Hess f
\approx \mathbf g$ and $\inf_{x\in V} |\nabla f(x)|>0$ locally}
\label{introducao}

${}$
\smallskip

Throughout this subsection we are going to denote the
Laplacian-Beltrami operator $\Delta_{\mathcal{M}}$ by $\Delta$ and
the tangential-gradient $\nabla_T$ by $\nabla$. Let $\mathcal M$
be a compact $n$-dimensional Riemannian manifold (without
boundary) with Riemmanian metric $\mathbf g$ of class $C^2$. Let
$\nabla$ denote the Levi-Civita connection. Fix $p \in
\mathcal{M}$. Our aim is to construct a function $f:V_p
\rightarrow \mathbb R$ such that $Hess f \approx \mathbf g$ and
$\underset{x\in V_p}\inf |\nabla f(x)|>0$, where $V_p$ is a
neighborhood of $p$ and the Hessian of $f$ is seen as a bilinear
form defined on the tangent space $T_p\mathcal M$ of $\mathcal M$
at $p$.

We begin with an orthonormal basis $(e_1,\ldots, e_n)$ of $T_pM$.
Define a normal coordinate system $(x_1,\ldots,x_n)$ in a
neighborhood $\widetilde V_p$ of $p$ such that $\partial/\partial
x_i(p)=e_i(p)$ for every $i=1,\ldots,n$. It is well known that in
this coordinate system we have that $\Gamma^k_{ij}(p)=0$, where
$\Gamma^k_{ij}$ are the Christoffel symbols with respect to
$(x_1,\ldots,x_n)$ (See, for instance, \cite{doCarmo}).

The Hessian with respect to $(x_1,\ldots,x_n)$ is given by
\[
Hess f \left( \frac{\partial}{\partial x_i},\frac{\partial}{\partial
x_j} \right)=\frac{\partial^2 f}{\partial x_i \partial
x_j}-\sum_{k=1}^n \Gamma^k_{ij} \frac{\partial f}{\partial x_k}.
\]

The Laplacian of $f$ is the trace of the Hessian with respect to the
metric $\mathbf g$. If $\mathbf g_{ij}$ denote the components of the
Riemannian metric with respect to $(x_1,\ldots,x_n)$ and $\mathbf
g^{ij}$ are the components of the inverse matrix of $\mathbf
g_{ij}$, then the Laplacian of $f$ is given by
\[
\Delta f=\sum_{i,j} \mathbf g^{ij} Hess f\left(
\frac{\partial}{\partial x_i},\frac{\partial}{\partial x_j}\right).
\]

Consider the function $f:\widetilde V_p \rightarrow \mathbb R$
defined by
\[
f(x)=x_1+\frac{1}{2}\sum_{i=1}^n x_i^2.
\]
It is immediate that $\Delta f(p)=n$ and $|\nabla f(p)|=1$.
Moreover, $Hess f(p)=\mathbf g(p)$, which implies that
\[
Hess f(p)(v,v)=| v |^2_p.
\]

We are interested in finding a neighborhood $V_p
 \subset \widetilde V_p$ of $p$ and a
strictly positive constant $C$ such that
\[
C \int_0^T \int_{V_p} \left(| \nabla u |^2 + u_t^2\right)
d\mathcal M dt
\]
\begin{equation}
\label{equacao1} \leq \int_0^T \int_{V_p} \left[Hess f(\nabla
u,\nabla u) + \left(\alpha-\frac{\Delta f}{2}\right) | \nabla
u|^2+ \left(\frac{\Delta f}{2} - \alpha \right)u_t^2 \right]
d\mathcal M\,dt,
\end{equation}
for some $\alpha \in \mathbb R$. We claim that if we consider
$\alpha=\frac{n}{2}-\frac{1}{2}$ and $C=1/4$ we obtain the desired
inequality, what means that it is enough to prove that there exist
$V_p\subset \widetilde V_p$ verifying
\begin{equation}
\label{desigualdadegradiente} \int_0^T \int_{V_p} Hess f(\nabla
u,\nabla u) + \left(\frac{n}{2}-\frac{3}{4}-\frac{\Delta
f}{2}\right) | \nabla u|^2 d\mathcal M dt\geq 0
\end{equation}
and
\begin{equation}
\label{desigualdadeescalar} \int_0^T \int_{V_p} \left(\frac{\Delta
f}{2} - \frac{n}{2}+\frac{1}{4} \right)u_t^2 d\mathcal M dt\geq 0.
\end{equation}

In order to prove the existence of a subset $V_p\subset \widetilde
V_p$ where (\ref{desigualdadegradiente}) holds, let $\theta_1$ be
the smooth field of symmetric bilinear form on $\widetilde V_p$
defined as
\[
\theta_1(X,Y)=Hess f(X,Y)+\left(\frac{n}{2}-\frac{3}{4}-\frac{\Delta
f}{2}\right)\mathbf g(X,Y)
\]
where $X$ and $Y$ are vector fields on $\widetilde V_p$. It is
clearly a positive definite bilinear form on $p$ since $Hess
f(p)(X,Y)=\mathbf g(p)(X,Y)$ and
\[
\theta_1(p)(X,Y)=\frac{1}{4} \mathbf g(p)(X,Y).
\]

Therefore, there exist a neighborhood $\widehat V_p$ such that
$\theta_1$ is positive definite and
\[
\int_0^T \int_{\widehat V_p} Hess f(\nabla u,\nabla u) +
\left(\frac{n}{2}-\frac{3}{4}-\frac{\Delta f}{2}\right) | \nabla
u|^2 d\mathcal M dt\geq 0.
\]

To prove the existence of $\overset{\smallsmile} V_p\subset
\widetilde V_p$ such that (\ref{desigualdadeescalar}) holds is
easier. It is enough to notice that at $p$ we have that
\[
\left(\frac{\Delta f(p)}{2} - \frac{n}{2}+\frac{1}{4}
\right)=\frac{1}{4}
\]
and the existence of $\overset{\smallsmile} V_p\subset \widetilde
V_p$ such that (\ref{desigualdadeescalar}) holds is immediate.
Furthermore we can eventually choose a smaller $V_p$ such that
$\underset{x\in V_p} \inf |\nabla f(x)|>0$. Therefore the
existence of $V_p \subset \widetilde V_p$ such that
$\underset{x\in V_p} \inf |\nabla f(x)|>0$ and (\ref{equacao1})
holds is settled.

\subsection{A function $f$ that satisfies Inequality (\ref{equacao1})
and $\inf_{x\in V} |\nabla f(x)|>0$ in a wide domain}
\label{hessianofuncaoquasemetrica}

${}$

\smallskip

 In what follows, $\bar V$ denotes the closure of
$V$ and $\partial V$ denotes the boundary of $V$. When $\bar V
\subset W$ is bounded, we say that $V$ is compactly contained in
$W$ and we denote this relationship by $V \subset \subset W$.

\begin{theorem}
\label{principal} Let $(\mathcal M,\mathbf g)$ be a two
dimensional Riemannian manifold. Then, for every $\epsilon
> 0$, there exist a finite family $\{V_i\}_{i=1\ldots k}$ of open
sets with smooth boundary, smooth functions $f_i:\bar V_i
\rightarrow \mathbb R$ and a constant $C>0$ such that
\begin{enumerate}
\item The subsets $\bar V_i$ are pairwise disjoint;
\item $\mathrm{vol}(\bigcup_{i=1}^k V_i)\geq
\mathrm{vol}(M)-\epsilon$;
\item Inequality (\ref{equacao1}) holds for every $f_i$;
\item $\underset{x\in V_i}\inf |\nabla f(x)|>0$ for every $i=1, \cdots,k$.
\end{enumerate}
\end{theorem}

\begin{proof}
First of all, it is possible to get open subsets $\{\widetilde
W_j\}_{j=1,\ldots,s}$ with smooth boundaries and a family of
smooth functions $\{\widetilde f_j:\widetilde W_j \rightarrow
\mathbb R\}_{j=1,\ldots,s}$ such that $\{\widetilde
W_j\}_{j=1,\ldots,s}$ is a cover of $\mathcal M$ and each
$\widetilde f_j$ satisfies Inequality (\ref{equacao1}). Moreover,
we can choose $\widetilde W_j$ in such a way that their boundaries
intercept themselves transversally and three or more boundaries do
not intercept themselves at the same point.

Set by $A:=\bigcup_{j=1}^s \partial \widetilde W_j$. Then,
$\mathcal{M}\backslash A$ is a disjoint union of connected open
sets $\bigcup_{i=1}^k W_i$ such that $\partial W_i$ is a piecewise
smooth curve.

Each $W_i$ is contained in some $\widetilde W_j$. Therefore, for
each $W_i$, choose a function $\hat f_i:=\widetilde f_j
\vert_{W_i}$.

The open subsets $V_i$, $i=1,\ldots,k$, we are looking for are
subsets of $W_i$. We can choose them in such a way that
\begin{enumerate}
\item $V_i \subset \subset W_i$;
\item $\partial V_i$ is smooth;
\item $\mathrm{vol}(W_i)-\mathrm{vol}(V_i)<\epsilon/k$.
\end{enumerate}

Finally, if we set $f_i=\hat f_i\vert_{\bar V_i}$, we prove the
theorem.

\end{proof}

\begin{theorem}
Let $(\mathcal M,\mathbf g)$ be a two-dimensional Riemannian
manifold. Fix $\epsilon
>0$. Then, there exist a smooth function $f:M \rightarrow \mathbb R$
such that inequality (\ref{equacao1}) and the condition
$\underset{x\in V}\inf |\nabla f(x)|>0$ hold in a subset $V$ with
$\mathrm{vol}(V)\geq \mathrm{vol}(\mathcal M)-\epsilon$.
\end{theorem}

\begin{proof}
Consider Theorem \ref{principal} and the constructions made in its
proof. Denote $\lambda:=\min\limits_{i\neq
j}\mathrm{dist}(V_i,V_j)>0$. Consider a tubular neighborhood
$V^\delta$ of $V=\cup_{i=1}^k V_{i}$ of the points whose distance
is less than or equal to $\delta<\lambda/4$. Then, it is possible
to define a smooth (cut-off) function $\eta:\mathcal M \rightarrow
\mathbb R$ such that
\[
\begin{array}{ccc}
\eta(x) & = & \left\{
\begin{array}{l}
1 \; \; \mathrm{if} \; \; x \in V \\
0 \; \; \mathrm{if} \; \; x \in \mathcal M\backslash V^\delta \\
\mathrm{between}\; \; 0 \;\; \mathrm{and} \;\; 1 \;\;
\mathrm{otherwise}.
\end{array}
\right.
\end{array}
\]

Now, notice that $f:\mathcal M\rightarrow \mathbb R$ defined by
\[
\begin{array}{ccc}
f(x) & = & \left\{
\begin{array}{l}
\hat f_i(x)\eta(x) \; \; \mathrm{if} \; \; x \in W_i; \\
0 \; \; \mathrm{otherwise}
\end{array}
\right.
\end{array}
\]
is smooth and satisfy inequality (\ref{equacao1}) and the
condition $\underset{x\in V}\inf |\nabla f(x)|>0$. In addition,
the inequality $\mathrm{vol}(V) \geq \mathrm{vol}(\mathcal
M)-\epsilon$ holds, which settles the theorem.
\end{proof}

%%%%%%%%%%%%%%%%%%%%%%%%%%%%%%%%%%%%%%%%%%%%%%%%%%%%%%%%%%%%%%%%%%%%%%%%%%%%%%%%%%%%%%%%%%%%%%%%%%%%%%%%

 We  denote
\begin{eqnarray}\label{Chi}
\chi=\left[ \int_{\mathcal{M}}u_{t}\, \nabla_{T}f\cdot \nabla_T
u\,d\mathcal{M}\right]_{0}^{T}+\alpha\left[ \int_{\mathcal{M}
}u_{t}\, u\,d\mathcal{M}\right] _{0}^{T}.
\end{eqnarray}

Next we will estimate some terms in (\ref{A}). Let us define
\begin{eqnarray}\label{R}
R:= \max_{x \in \mathcal{M}} |\nabla_T f(x)|.
\end{eqnarray}

\smallskip

\noindent{\em Estimate for $I_1:= \int_{0}^{T}\int_{\mathcal{M}
}a(x)\,g(u_{t})(\nabla_{T} f \cdot\nabla_T u)d\mathcal{M} dt.$}

\smallskip

By Cauchy-Schwarz inequality, taking (\ref{R}) into account and
considering the inequality $ab\leq \frac{a^2}{4\zeta} + \zeta
b^2$, where $\zeta$ is a positive number, we obtain
\begin{eqnarray}\label{5.2.14}
|I_1| \leq \frac{||a||_{L^\infty(\mathcal{M})}R^2}{\zeta} \int_0^T
\int_{\mathcal{M}} a(x)|g(u_t)|^2d\mathcal{M} dt + 2\zeta \int_0^T
E(t)\,dt.
\end{eqnarray}

\noindent{\em Estimate for $I_2= \alpha\int_0^T
\int_{\mathcal{M}}a(x)\,g(u_t)u\, d\mathcal{M} dt.$}

Similarly we infer
\begin{eqnarray}\label{5.2.17}
|I_2| \leq
\frac{||a||_{L^\infty(\mathcal{M})}\alpha^2\lambda_1^{-1}}{16\zeta}\int_0^T
\int_{\mathcal{M}} a(x)|g(u_t)|^2\,d\mathcal{M} dt +  2\zeta
\int_0^T E(t)\,dt,
\end{eqnarray}
where $\lambda_1$ comes from the Poincaré inequality given in
(\ref{Poincare}).

\smallskip

Choosing $\zeta$ sufficiently small and inserting (\ref{Chi}),
(\ref{5.2.14}) and (\ref{5.2.17}) into (\ref{A}) yields
\begin{eqnarray}\label{5.2.18}
\int_0^T E(t)\,dt &\leq& |\chi| +
C_1\int_0^T\int_{\mathcal{M}}a(x)\,(g(u_t))^2
d\mathcal{M} dt\\
&+&C_1\int_0^T\int_{\mathcal{M}\backslash V}[ |\nabla_T u|^2 +
a(x)\,u_t^2]\,d\mathcal{M} dt\nonumber
\end{eqnarray}
where
\begin{eqnarray*}
C_1:= C_1\left\{C,||a||_{L^\infty(\mathcal{M})},\lambda_1^{-1},R,
a_0^{-1},n\right\}.
\end{eqnarray*}

It remains to estimate the quantity
$\int_0^T\int_{\mathcal{M}\backslash V}|\nabla_T
u|^2\,d\mathcal{M} dt$ in terms of the damping term
$\int_0^T\int_{\mathcal{M}}[a(x)\,|g(u_t)|^2 + a(x)\,|u_t|^2
]\,d\mathcal{M} dt$. For this purpose we have to build a
``cut-off'' function $\eta_\varepsilon$ on a specific neighborhood
of $\mathcal M \backslash V$. First of all, define $\tilde \eta:
\mathbb R \rightarrow \mathbb R$ such that
\[
\begin{array}{ccc}
\tilde\eta(x) & = & \left\{
\begin{array}{clc}
1 & \mathrm{if} & x \leq 0 \\
(x-1)^2 & \mathrm{if} & x \in [1/2,1] \\
0 & \mathrm{if} & x>1
\end{array}
\right.
\end{array}
\]
and it is defined on $(0,1/2)$ in such a way that $\tilde \eta$ is a
non-decreasing function of class $C^1$. For $\varepsilon > 0$, set
$\tilde \eta_\varepsilon(x):=\tilde\eta(x / \varepsilon)$. It is
straightforward that there exist a constant $M$ which does not
depend on $\varepsilon$ such that
\[
\frac{\vert \tilde\eta_\varepsilon^\prime(x)
\vert^2}{\tilde\eta_\varepsilon(x)} \leq \frac{M}{\varepsilon^2}
\]
for every $x < \varepsilon$.

Let $\mathcal M_*\supset \mathcal M\backslash V$ be an open subset
of $\mathcal M$ and let $\varepsilon
> 0$ be such that
\[
\tilde \omega_{\varepsilon}:=\{x \in \mathcal
M;\mathrm{dist}(x,\partial V)< \varepsilon\}
\]
is a tubular neighborhood of $\partial V$ and
$\omega_{\varepsilon}:= \tilde \omega_{\varepsilon} \cup
\mathcal{M}\backslash V $ is contained in $\mathcal M_{*}$. Define
$\eta_\varepsilon:\mathcal M \rightarrow \mathbb R$ as
\[
\begin{array}{ccc}
\eta_\varepsilon(x) & = & \left\{
\begin{array}{clc}
1 & \mathrm{if} & x \in \mathcal{M}\backslash V \\
\tilde\eta_\varepsilon(d(x,\mathcal{M}\backslash V)) & \mathrm{if}
& x \in \omega_{\varepsilon} \backslash
(\mathcal{M}\backslash V) \\
0 & \mathrm{otherwise}. &
\end{array}
\right.
\end{array}
\]

It is straightforward that $\eta_\varepsilon$ is a function of
class $C^1$ on $\mathcal M$ due to the smoothness of $\partial
(\mathcal{M}\backslash V)$ and $\partial \omega_\varepsilon$.
Notice also that
\begin{eqnarray}\label{eq:3.38'}
\frac{\vert \nabla_T \eta_\varepsilon
(x)\vert^2}{\eta_\varepsilon(x)}=\frac{\vert \tilde
\eta^\prime_\varepsilon (d(x,\mathcal M_2))\vert^2}{\tilde
\eta_\varepsilon(d(x,\mathcal M_2))} \leq \frac{M}{\varepsilon^2}
\end{eqnarray}
for every $x\in \omega_\varepsilon$. In particular, $\frac{\vert
\nabla_T \eta_\varepsilon \vert^2}{\eta_\varepsilon} \in
L^\infty(\omega_\varepsilon)$.

Taking $\xi=\eta_\varepsilon$ in the identity (\ref{identity II}) we
obtain
\begin{eqnarray} \label{5.3.1}
&&\int_0^T \int_{\omega_\varepsilon} \eta_\varepsilon |\nabla_T
u|^2 d\mathcal{M} dt \\
&=& -\left[ \int_{\omega_\varepsilon} u_t u \eta_\varepsilon
\,d\mathcal{M}\right]_0^T
+\int_0^T \int_{\omega_\varepsilon} \eta_\varepsilon |u_t|^2\,d\mathcal{M} \nonumber\\
&-& \int_0^T \int_{\omega_\varepsilon}u(\nabla_T u \cdot \nabla_T
\eta_\varepsilon) \,d\mathcal{M} dt - \int_0^T
\int_{\omega_\varepsilon} a(x)\, g(u_t)u \eta_\varepsilon\,
d\mathcal{M} dt.\nonumber
\end{eqnarray}

Next we will estimate the terms on the RHS of (\ref{5.3.1}).

\smallskip
\noindent{\em Estimate for $K_1:= \int_0^T \int_{\omega_\varepsilon}
\eta_\varepsilon |u_t|^2\,d\mathcal{M} dt$}
\smallskip

From (\ref{eq:2.2}), since $\eta_\varepsilon \leq 1$ and
$\omega_{\varepsilon} \subset \mathcal{M}_{\ast}$, where the damping
lies, we deduce
\begin{eqnarray} \label{5.3.4}
K_1 \leq a_0^{-1} \int_0^T \int_{\mathcal{M}} a(x)\,
u_t^2\,d\mathcal{M}\,dt.
\end{eqnarray}

\smallskip
\noindent{\em Estimate for $K_2:= - \int_0^T
\int_{\omega_{\varepsilon}} a(x)\,g(u_t)u \eta_\varepsilon
\,d\mathcal{M} dt.$}
\smallskip

Taking into account the Cauchy-Schwarz inequality, the inequality
$ab \leq \frac{1}{4\alpha}a^2 + \alpha b^2$ and (\ref{Poincare})
we obtain
\begin{eqnarray} \label{5.3.5}
|K_2|\leq
\frac{\lambda_1^{-1}||a||_{L^{\infty}(\mathcal{M})}}{4\alpha}\int_0^T\int_{\mathcal{M}}a(x)\,|g(u_t)|^2\,d\mathcal{M}
+ 2\alpha \int_0^T E(t)\,dt,
\end{eqnarray}
where $\alpha$ is a positive constant.

\smallskip
\noindent{\em Estimate for $K_3:=  \int_0^T
\int_{\omega_{\varepsilon}}u(\nabla_T u \cdot \nabla_T
\eta_\varepsilon) d\mathcal{M} dt.$}
\smallskip

Considering (\ref{eq:3.38'}) and applying Cauchy-Schwarz inequality,
we can write
\begin{eqnarray}\label{5.3.6}
|K_3|&\leq& \frac12 \int_0^T \left[
\int_{\omega_{\varepsilon}}\eta_\varepsilon |\nabla_T
u|^2\,d\mathcal{M} + \int_{\omega_{\varepsilon}} \frac{|\nabla_T
\eta_\varepsilon |^2}{\eta_\varepsilon}|u|^2\,d\mathcal{M}\right]dt\\
&\leq&\frac12 \int_0^T \left[
\int_{\omega_{\varepsilon}}\eta_\varepsilon |\nabla_T
u|^2\,d\mathcal{M} +
\frac{M}{\varepsilon^2}\int_{\omega_{\varepsilon}}
|u|^2\,d\mathcal{M}\right]dt. \nonumber
\end{eqnarray}

Combining (\ref{5.3.1})-(\ref{5.3.6}) we arrive to the following
inequality
\begin{eqnarray}\label{5.3.7}\qquad
\frac12 \int_0^T \int_{\omega_{\varepsilon}} \eta_\varepsilon
|\nabla_T u|^2 \,d\mathcal{M} dt &\leq& |\mathcal{Y}| +
\frac{\lambda_1^{-1}||a||_{L^{\infty}(\mathcal{M})}}{4\alpha}\int_0^T\int_{\mathcal{M}}a(x)\,|g(u_t)|^2\,d\mathcal{M}\\
&+& 2\alpha \int_0^T E(t)\,dt+ \frac{M}{2\varepsilon^2}\int_0^T
\int_{\omega_{\varepsilon}} |u|^2\,d\mathcal{M}\,dt
,\nonumber\\
&+& a_0^{-1} \int_0^T \int_{\mathcal{M}} a(x)\,
u_t^2\,d\mathcal{M}\,dt. \nonumber
\end{eqnarray}
where
\begin{eqnarray}\label{Z''}
\mathcal{Y}:=-\left[ \int_{\omega_{\varepsilon}} u_t u
\eta_\varepsilon \,d\mathcal{M}\right]_0^T.
\end{eqnarray}

Thus, combining (\ref{5.3.7}) and (\ref{5.2.18}), have in mind
that
\begin{eqnarray*}
\frac12 \int_0^T \int_{\mathcal{M}\backslash V} |\nabla_T u|^2
\,d\mathcal{M} dt \leq \frac12 \int_0^T
\int_{\omega_{\varepsilon}} \eta_\varepsilon |\nabla_T u|^2
\,d\mathcal{M} dt
\end{eqnarray*}
and choosing $\alpha$ small enough, we deduce
\begin{eqnarray}\label{5.3.8}
&&\int_0^T E(t)\,dt \leq  |\chi|  +  C_1|\mathcal{Y}|\qquad \\
&&+C_2\int_0^T\int_{\mathcal{M}}[a(x)\,|g(u_t)|^2 + a(x)\,|u_t|^2
]\,d\mathcal{M} dt\nonumber\\
&&+\frac{M C_2}{\varepsilon^2}\int_0^T \int_{\omega_{\varepsilon}}
\,|u|^2\,d\mathcal{M}\,dt,\nonumber
\end{eqnarray}
where $C_2=C_2(C_1,\lambda_1^{-1},||a||_{L^\infty(\mathcal
{M})},a_0^{-1})$.

On the other hand, from (\ref{Chi}),  (\ref{Z''}) and (\ref{3.12})
the following estimate holds
\begin{eqnarray}\label{5.3.9}
|\chi|+  2C_2 |\mathcal{Y}|  &\leq&
C(E(0)+E(T))\\
&=& C\left[ 2\,E(T) + \int_0^T \int_{\mathcal{M}}
a(x)\,g(u_t)\,u_t\,d\mathcal{M}\right],\nonumber
\end{eqnarray}
where $C$ is a positive constant which depends on $R$.

Then, (\ref{5.3.8}) and (\ref{5.3.9}) yield
\begin{eqnarray}\label{5.3.10}
T\,E(T)&\leq& \int_0^T E(t)\,dt\\
&\leq& C\,E(T) + C
\left[\int_0^T\int_{\mathcal{M}}[a(x)\,|g(u_t)|^2 +
a(x)\,|u_t|^2 ]\,d\mathcal{M} dt\right]\nonumber\\
&+& C\int_0^T \int_{\omega_{\varepsilon}}
\,|u|^2\,d\mathcal{M}\,dt,\nonumber
\end{eqnarray}
where $C$ is a positive constant which depends on $a_0, \lambda_1,
 R, ||a||_{L^\infty(\mathcal
{M})}, n $ and $ \frac{M}{\varepsilon^2}$.

%%%%%%%%%%%%%%%%%inicio do lema %%%%%%%%%%%%%%%%%%%%%%%%%%%%%%%%
Our aim is to absorb the last term on the RHS of (\ref{5.3.10}).
In order to do this, let us consider the following lemma, where
$T_0$ is a positive constant which is sufficiently large for our
purpose.

\medskip
\begin{lemma}\label{lemma3.5}
Under the hypothesis of Theorem \ref{Theo. 3.1} and for all
$T>T_0$, there exists a positive constant $C(T_0,E(0))$ such that
if $(u,u_t)$ is the solution of (\ref{1.1}) with weak initial
data, we have
\begin{eqnarray}\label{3.53}\quad
\int_0^T \int_{\mathcal{M}} |u|^2 \,d\mathcal{M}\,dt\leq
C(T_0,E(0))
\left\{\int_0^T\int_{\mathcal{M}}\left(a(x)\,g^2(u_t))+ a(x)
u^2_t\right)d\mathcal{M}\,dt \right\}.
\end{eqnarray}
\end{lemma}

{\em Proof:} We argue by contradiction exactly as in Lasiecka and
Tataru's work \cite{Lasiecka-Tataru}. For simplicity we shall
denote $u' := u_t$. Let us suppose that (\ref{3.53}) is not
verified and let $\{u_k(0),u'_k(0)\}$ be a sequence of initial
data where the corresponding solutions $\{u_k\}_{k\in \mathbb{N}}$
of (\ref{1.1}), with $E_k(0)$ assumed uniformly bounded in $k$,
verifies
\begin{eqnarray}\label{3.54}
\lim_{k \rightarrow +\infty}\frac{\int_0^T
||u_k(t)||_{L^2(\mathcal{M})}^2
dt}{\int_0^T\int_{\mathcal{M}}\left( a(x)\,g^2(u_k')+
a(x)\,u_k'^2\right)d\mathcal{M}\,dt}=+\infty,
\end{eqnarray}
that is
\begin{eqnarray}\label{3.55}
\lim_{k \rightarrow +\infty}\frac{\int_0^T\int_{\mathcal{M}}\left(
a(x)\,g^2(u_k')+ a(x)\,u_k'^2\right)d\mathcal{M}\,dt}{\int_0^T
||u_k(t)||_{L^2(\mathcal{M})}^2 dt}=0.
\end{eqnarray}

Since $E_k(t) \leq E_k(0)\leq L$, where $L$ is a positive
constant, we obtain a subsequence, still denoted by $\{u_k\}$ from
now on, which verifies the convergence
\begin{eqnarray}
&& u_k \rightharpoonup u \hbox{ weakly in } H^1(\Sigma_T),\label{3.56}\\
&& u_k \rightharpoonup u \hbox{ weak star in } L^{\infty}(0,T;
V),\label{3.57}\\
&& u_k' \rightharpoonup u' \hbox{ weak star in } L^{\infty}(0,T;
L^2(\mathcal{M})).\label{3.58}
\end{eqnarray}

Employing compactness results we also deduce that
\begin{eqnarray}
u_k \rightarrow u \hbox{ strongly in }
L^2(0,T;L^2(\mathcal{M})).\label{3.60}
\end{eqnarray}

\medskip

At this point we will divide our proof into two cases, namely,
$u\ne 0$ and $u=0$.

\medskip

(i) Case (I): $u \ne 0.$

We also observe that from (\ref{3.55}) and (\ref{3.60}) we have
\begin{eqnarray}\label{3.62}
\lim_{k \rightarrow +\infty}\int_0^T\int_{\mathcal M}\left(
a(x)\,g^2(u_k')+ a(x)\,u_k'^2\right)d\mathcal{M}\,dt=0
\end{eqnarray}

Passing to the limit in the equation, when $k \rightarrow
+\infty$, we get,
\begin{equation}\label{3.63}
\left\{
\begin{aligned}
u_{tt}-\Delta_{\mathcal{M}} \,u&=& 0 \hbox{ on }~\mathcal{M} \times (0,T)\\
u_t &=& 0 \hbox{ on } \mathcal{M}_{\ast}\times(0,T),
\end{aligned}
\right.
\end{equation}
and for $u_t=v$, we obtain, in the distributional sense
\begin{equation*}
\left\{
\begin{aligned}
v_{tt}-\Delta_{\mathcal{M}} \,v&=& 0 \hbox{ on }~\mathcal{M} \times (0,T),\\
v &=& 0 \hbox{ on } \mathcal{M}_{\ast}\times (0,T).
\end{aligned}
\right.
\end{equation*}

From uniqueness results due to Triggiani and Yao \cite{Trigianni}
we conclude that $v\equiv 0$, that is, $u_t=0$. Indeed, let $V_i$
as in Theorem \ref{principal} and $\Gamma=\partial V_i$, which is
a smooth curve contained in $\mathcal M_\ast$. Since $v\equiv 0$
on $\mathcal{M}_{\ast}$ we deduce that $v=\partial_{\nu}v=0$ on
$\Gamma$. Employing Triggiani and Yao's uniqueness results to the
compact manifold $\bar V_i$ with boundary $\Gamma$ we infer that
$v\equiv 0$ on $\bar V_i$, for each $i=1, \cdots,k$.  Therefore,
$v \equiv0$ on $\mathcal{M}$ as we desired to prove. Returning to
(\ref{3.63}) we obtain the following elliptic equation  a.e. in $
(0,T)$ given by
\begin{equation*}
\left\{
\begin{aligned}
\Delta_{\mathcal{M}} \,u&=& 0 \hbox{ on }~\mathcal{M}\\
u_t &=& 0 \hbox{ on } \mathcal{M},
\end{aligned}
\right.
\end{equation*}
which implies that $u=0$, which is a contradiction.

\medskip
(ii) Case (II): $u = 0.$

\medskip

Defining
\begin{eqnarray}\label{3.66}
c_k := \left[ \int_0^T \int_{\mathcal{M}} |u_k|^2
d\mathcal{M}\,dt\right]^{1/2}
\end{eqnarray}
and
\begin{eqnarray}\label{3.67}
\overline{u}_k:= \frac{1}{c_k}\,u_k,
\end{eqnarray}
we obtain
\begin{eqnarray}\label{3.68}
&\int_0^T\int_{\mathcal{M}} |\overline{u}_k|^2 d\mathcal{M}\,dt =
\int_0^T\int_{\mathcal{M}} \frac{|u_k|^2}{c_k^2}d\mathcal{M}\,dt=
\frac{1}{c_k^2}\int_0^T \int_{\mathcal{M}} |u_k|^2
d\mathcal{M}\,dt=1. &
\end{eqnarray}

Setting
\begin{eqnarray*}
\overline{E}_k(t)&:=& \frac12 \int_{\mathcal{M}}
|\overline{u}'_k|^2\,d\mathcal{M} + \frac12 \int_{\mathcal{M}}
|\nabla \overline{u}_k|^2\,d\mathcal{M} ,
\end{eqnarray*}
we deduce that
\begin{eqnarray}\label{3.69}
\overline{E}_k(t)= \frac{E_k(t)}{c_k^2}.
\end{eqnarray}

Recalling (\ref{5.3.10}) we obtain, for  $T$ large enough, that
\begin{eqnarray*}
E(T) \leq \hat{C}\left[\int_0^T \int_{\mathcal{M}} (a(x)\,
g^2(u_t)+ a(x)\,u_t^2 )\,d\mathcal{M}\,dt +
\int_0^T\int_{\mathcal{M}} |u|^2\,d\mathcal{M}\,dt\right],
\end{eqnarray*}
and employing the identity $E(T)-E(0)= -\int_0^T\int_{\mathcal{M}}
a(x)\,g(u_t)\,u_t\,d\mathcal{M}\,dt$, we get
\begin{eqnarray*}
E(t) \leq E(0) \leq \tilde{C}\left[\int_0^T \int_{\mathcal{M}}
(a(x)\,g^2(u_t)+ a(x)\,u_t^2 )\,d\mathcal{M}\,dt +
\int_0^T\int_{\mathcal{M}} |u|^2\,d\mathcal{M}\,dt\right],
\end{eqnarray*}
for all $t\in (0,T)$, where $T$ is sufficiently large . The last
inequality and (\ref{3.69}) yield
\begin{eqnarray}\label{3.70}
\overline{E}_k(t):=\frac{E_k(t)}{c_k^2} \leq \tilde{C}
\left[\frac{\int_0^T \int_{\mathcal{M}}(a(x)\,g^2(u_k')+
a(x)\,u_k'^2)}{\int_0^T
\int_{\mathcal{M}}|u_k|^2\,d\mathcal{M}\,dt}+1 \right].
\end{eqnarray}

From (\ref{3.55}) and (\ref{3.70}) we conclude that there exists a
positive constant $\hat{M}$ such that
\begin{eqnarray*}
\overline{E}_k(t):=\frac{E_k(t)}{c_k^2} \leq \hat{M}, ~\hbox{ for
all } t\in[0,T] ~\hbox{ and for all } k\in \mathbb{N},
\end{eqnarray*}
that is,
\begin{eqnarray}\label{3.71}\qquad
\frac12 \int_{\mathcal{M}} |\overline{u}'_k|^2\,d\mathcal{M} +
\frac12 \int_{\Omega} |\nabla \overline{u}_k|^2\,d\mathcal{M}\leq
\hat{M}, ~\hbox{ for all } t\in[0,T] ~\hbox{ and for all } k\in
\mathbb{N}.
\end{eqnarray}

For a subsequence $\{\overline{u}_k\}$, we obtain
\begin{eqnarray}
&&\overline{u}_k \rightharpoonup \overline{u} \hbox{ weak star
in }L^{\infty}(0,T;V),\label{3.72}\\
&&\overline{u}'_k  \rightharpoonup \overline{u}' \hbox{ weak
star in } L^{\infty}(0,T;L^2(\mathcal{M})),\label{3.73}\\
&&\overline{u}_k  \rightarrow \overline{u} \hbox{ strongly in }
L^2(0,T;L^2(\mathcal{M})).\label{3.74}
\end{eqnarray}

We observe that from (\ref{3.55}) we deduce
\begin{eqnarray}\label{3.75}\qquad
\lim_{k \rightarrow +\infty}\int_0^T \int_{\mathcal{M}}
\frac{a(x)\,g^2(u_k')}{c_k^2}\,d\mathcal{M}\,dt=0~\hbox{ and }
\lim_{k \rightarrow +\infty}\int_0^T \int_{\mathcal{M}}
a(x)\,|\overline{u}_k'|^2\,d\mathcal{M}\,dt=0.
\end{eqnarray}

In addition, $\overline{u}_k$ satisfies the equation
\begin{equation*}
\left.
\begin{aligned}
\overline{u}_k'' - \Delta_{\mathcal{M}}\overline{u}_k +
a(x)\,\frac{g(u_k')}{c_k}&=&0\quad \hbox{ on } \mathcal{M} \times
(0,T).
\end{aligned}
\right.
\end{equation*}

Passing to the limit when $k \rightarrow +\infty$ and taking the
above convergences into account, we obtain
\begin{equation}\label{P2}
\left\{
\begin{aligned}
\overline{u}'' - \Delta_{\mathcal{M}}\overline{u} &=&0\quad
\hbox{ on } \mathcal{M} \times (0,T),\\
\overline{u}'&=&0\quad \hbox{ on }\mathcal{M}_{\ast} \times (0,T).
\end{aligned}
\right.
\end{equation}

Then, $v= \overline{u}_t$ verifies, in the distributional sense
\begin{equation*}
\left\{
\begin{aligned}
v_{tt}-\Delta_{\mathcal{M}} \,v&=& 0 \hbox{ on }~\mathcal{M}\\
v &=& 0 \hbox{ on } \mathcal{M}_{\ast}.
\end{aligned}
\right.
\end{equation*}

Applying, again, uniqueness results due to Triggiani and Yao
\cite{Trigianni}, it results that $v=\overline{u}_t=0$. Returning
to (\ref{P2}) we have a.e. in $(0,T)$ that
\begin{equation*}
\left\{
\begin{aligned}
\Delta_{\mathcal{M}} \,\overline{u}&=& 0 \hbox{ on }~\mathcal{M}\\
\overline{u}_t &=& 0 \hbox{ on } \mathcal{M}.
\end{aligned}
\right.
\end{equation*}

We deduce that $\overline{u}=0$, which is a contradiction in view
of (\ref{3.68}) and (\ref{3.74}). The lemma is settled. \quad
$\square$
%%%%%%%%%%%%%%%%%%fim do lema %%%%%%%%%%%%%%%%%%%%%%%%%%%%%%%

\smallskip

Inequalities (\ref{5.3.10}) and (\ref{3.53}) lead us to the
following result.
\smallskip

\noindent \textbf{Proposition 5.2.2:} \textit{For }$T>0$\textit{\
large enough, the solution }$\left[ u,u_{t}\right] $\textit{\ of
}(\ref{3.1}) \textit{satisfies}
\begin{equation}
E(T)\leq C\,\int_{0}^{T}\int_{\mathcal{M}}\left[ a(x)\left\vert
u_{t}\right\vert ^{2}+ a(x)\left\vert g\left( u_{t}\right)
\right\vert ^{2}\right] d\mathcal{M} dt \label{5.3.18}
\end{equation}%
\textit{where the constant }$C=C(T_0, E(0), C, a_0, \lambda_1, R,
||a||_{L^{\infty}(\mathcal{M})},n, \frac{M}{\varepsilon^2}).$

\subsection{Conclusion of Theorem 3.1}

In what follows we will proceed exactly as in Lasiecka and
Tataru's  work\cite{Lasiecka-Tataru}(see Lemma 3.2 and Lemma 3.3
of the referred paper) adapted to our context. Let $\Sigma:=
\mathcal{M}\times (0,T)$,
\begin{eqnarray*}
\Sigma _{\alpha } &=&\left\{ \left( t,x\right) \in \Sigma /\text{ \ }%
\left\vert u_{t}\right\vert >1\text{ \ a. e.}\right\} , \\
\Sigma _{\beta } &=&\Sigma \backslash \Sigma _{\alpha }.
\end{eqnarray*}

Then using hypothesis $(iii)$ in Assumption \ref{ass:1}, we obtain%
\begin{equation}
\int_{\Sigma _{\alpha }}a(x)\left( \left[ g\left( u_{t}\right)
\right] ^{2}+\left( u_{t}\right) ^{2}\right) d\Sigma _{\alpha
}\leq \left( k^{-1}+K\right) \int_{\Sigma _{\alpha }}a(x)g\left(
u_{t}\right) u_{t}d\Sigma _{\alpha }.  \label{5.4.1}
\end{equation}%

Moreover, from (\ref{4.1'})
\begin{equation}
\int_{\Sigma _{\beta }}a(x)\left( \left[ g\left( u_{t}\right)
\right] ^{2}+\left( u_{t}\right) ^{2}\right) d\Sigma _{\beta }\leq
(1+||a||_{\infty})\int_{\Sigma _{\beta }}h\left(a(x) g\left(
u_{t}\right) u_{t}\right) d\Sigma _{\beta }. \label{5.4.2}
\end{equation}%

Then, by Jensen's inequality
\begin{eqnarray}
(1+||a||_{\infty})\int_{\Sigma _{\beta }}h\left( g\left(
u_{t}\right) u_{t}\right) d\Sigma _{\beta } &\leq
&(1+||a||_{\infty})meas\left( \Sigma\right) h\left(
\frac{1}{meas\left( \Sigma\right) }\int_{\Sigma}a(x)g\left(
u_{t}\right) u_{t}d\Sigma
\right) \smallskip  \notag \\
&=&(1+||a||_{\infty})meas\left( \Sigma\right) r\left(
\int_{\Sigma}a(x)g\left( u_{t}\right) u_{t}d\Sigma\right) ,
\label{5.4.3}
\end{eqnarray}
where $r\left( s\right) =h\left( \frac{s}{meas\left( \Sigma
\right) } \right)$ is defined in (\ref{4.2'}). Thus
\begin{eqnarray}
\int_{\Sigma}a(x)\left( \left[ g\left( u_{t}\right) \right]
^{2}+\left( u_{t}\right) ^{2}\right) d\Sigma &\leq &\left(
k^{-1}+K\right)
\int_{\Sigma}a(x) g\left( u_{t}\right)u_t\,d\Sigma \smallskip  \notag \\
&&+(1+||a||_{\infty})meas\left( \Sigma\right) r\left( \int_{\Sigma
}a(x)g\left( u_{t}\right) u_{t}\,d\Sigma\right) . \label{5.4.4}
\end{eqnarray}

Splicing, together, (\ref{5.3.18}) and (\ref{5.4.4}), we have
\begin{eqnarray}
E(T) &\leq &(1+||a||_{\infty})C\,\left[
\frac{K_{0}}{(1+||a||_{\infty})}\int_{\Sigma}a(x) g\left(
u_{t}\right)
u_{t}d\Sigma \right. \smallskip  \notag \\
&&\left. + meas\left( \Sigma\right) r\left(
\int_{\Sigma}a(x)\,g\left( u_{t}\right) u_{t}d\Sigma\right)
\right] , \label{5.4.5}
\end{eqnarray}
where $K_{0}=k^{-1}+K$. Setting
\begin{eqnarray*}
L &=&\frac{1}{C\,meas\left( \Sigma\right) (1+||a||_{\infty})},\smallskip \\
c &=&\frac{K_{0}}{meas\left( \Sigma \right)(1+||a||_{\infty}) },
\end{eqnarray*}
we obtain
\begin{eqnarray}
p\left[ E(T)\right] &\leq &\int_{\Sigma}a(x)\,g\left( u_{t}\right)
u_{t}\,d\Sigma  =E(0)-E(T),  \label{5.4.6}
\end{eqnarray}%
where the function $p$ is as defined in (\ref{4.3'}). To finish
the proof of Theorem 3.1, we invoke the following result from I.
Lasiecka et al. \cite{Lasiecka-Tataru}:

\smallskip

\noindent \textbf{Lemma B}:\textit{\ Let }$p$\textit{\ be a positive,
increasing function such that }$p(0)=0$\textit{. Since }$p$\textit{\ is
increasing we can define an increasing function }$q,$ $q(x)=x-(I+p)^{-1}%
\left( x\right) .$\textit{\ Consider a sequence }$s_{n}$\textit{\
of positive numbers which satisfies}
\begin{equation*}
s_{m+1}+p(s_{m+1})\leq s_{m}.
\end{equation*}
\indent\textit{Then }$s_{m}\leq S(m)$\textit{, where
}$S(t)$\textit{\ is a solution
of the differential equation}%
\begin{equation*}
\frac{d}{dt}S(t)+q(S(t))=0,\text{ \ }S(0)=s_{0}.
\end{equation*}%
\indent\textit{Moreover, if }$p(x)>0$\textit{\ for }$x>0$\textit{,
then }$\underset{ t\rightarrow \infty }{\lim }$\textit{\
}$S(t)=0.$

Taking into account the above result, we replace $T$ (resp. $0$)
in (\ref{5.4.6}) with $ m(T+1)$ (resp. \ $mT$) in order to get
\begin{equation}
E(m(T+1))+p\left( E(m(T+1))\right) \leq E(mT),\text{ \ for \
}m=0,1,.... \label{5.4.7}
\end{equation}

Applying Lemma B with $s_{m}=E(mT)$ results in
\begin{equation}
E(mT)\leq S(m),\text{ \ \ }m=0,1,....  \label{5.4.8}
\end{equation}%

Finally, using the inherent dissipativity of $E(t)$ given in
relation (\ref{3.12}), we have for $t=mT+\tau ,$ $0\leq \tau \leq
T,$
\begin{equation}
E(t)\leq E(mT)\leq S(m)\leq S\left( \frac{t-\tau }{T}\right) \leq
S\left( \frac{t}{T}-1\right) \text{ \ \ for \ }t>T\text{,}
\label{5.4.9}
\end{equation}
where we have used the fact that $S(.)$ is dissipative. The proof
of Theorem 3.1 is now completed.

\medskip
\vskip5pt {\bf Acknowledgements.} The authors would like to thank
Professor Roberto Triggiani  for his kind attention and helpful
remarks during the period when this paper was written.

\end{document}